\documentclass{article}
\usepackage[utf8]{inputenc}
\usepackage{lineno}
\usepackage{amssymb,mathrsfs, amsmath}
\usepackage{amsmath,amsthm}
\usepackage{amssymb}
\usepackage{latexsym}
\usepackage{amsfonts}
\usepackage{mathrsfs}
\usepackage{epsfig}
\usepackage{hyperref}
\usepackage{graphicx}
\usepackage{graphics}
\usepackage{pst-node}
\usepackage{stmaryrd}
\usepackage{cancel}
\usepackage {appendix}
\usepackage{titlefoot}
\usepackage{tikz-cd}
\newtheorem{thm}{Theorem}[section]

\newtheorem{defn}{Definition}[section]
\newtheorem{prop}{Proposition}[section]
\newtheorem{lem}{Lemma}[section]

\newtheorem{cor}{Corollary}[section]

\newtheorem{eg}{Example}[section]

\title{Cohomology and deformation of compatible Hom-Leibniz algebras}
\author{ Rinkila Bhutia, RB Yadav and Namita Behera  }
\date{}

\begin{document}

\maketitle
\begin{abstract} 
 In this article, we define  a suitable  graded Lie algebra whose Maurer-Cartan elements characterize the structure of compatible Hom-Leibniz algebras. We use this to define  cohomology of compatible Hom-Leibniz algebra. Using this cohomology we  study  infinitesimal deformations. Finally, we study  Nijenhuis operator and their relation with  compatible Hom-Leibniz algebras.
\end{abstract}

\vspace{0.3cm}

\vspace{0.3cm}


\unmarkedfntext{\newline\hspace*{2em} \textit{\thanks{Department of Mathematics, Sikkim University,
Sikkim-737102, INDIA,} ({\tt  rbhutia@cus.ac.in, rbyadav15@gmail.com, nbehera@cus.ac.in})}}

\section{Introduction}
Leibniz algebra is a non-commutative generalisation of Lie algebra. It was introduced and called D-algebra in papers by A. M. Bloch published in the 1960s to signify its relation with derivations. Later in 1993 J. L. Loday \cite{Loday} introduced the same structure and called it Leibniz algebra. The cohomology theory of Leibniz algebra with coefficients in a bimodule has been studied in \cite{LodayPiras}.
 Hom-Lie algebras were introduced by Hartwig, Larsson, and Silverstrov \cite{Hom}. Makhlouf and Silverstrov \cite{Makh} introduced the notion of Hom-Leibniz algebra, generalising Hom-Lie algebras. Hom-algebra structures have been widely studied since then.\\

Algebraic deformation theory was introduced by Gerstenhaber for rings and algebra in a series of papers \cite{MG1}-\cite{MG5}. Subsequently, algebraic deformation theory has been studied for different algebras by many authors. 
In \cite{Bala}, D. Balavoine studies the formal deformation of algebras using the theory of Maurer-Cartan elements in a graded Lie algebra. In particular, this approach is used to study the deformation of Leibniz algebra.\\ 

 Recently, cohomology and infinitesimal deformations of compatible Lie algebra and compatible associative algebra have been studied in \cite{CLieA} and \cite{CAA} respectively. Motivated by these works, we give a characterization of compatible Hom-Leibniz algebraand introduce  cohomology theory of such algebras.  We  define a compatible Hom-Leibniz algebra to be a pair of Hom-Leibniz algebras such that the linear combination of their algebraic structures is also a Hom-Leibniz algebra structure. 

 We define a graded Lie algebra whose Maurer-Cartan elements characterize the structure of compatible Hom-Leibniz algebras. We then study the cohomology of a compatible Hom-Leibniz algebra with coefficients in itself. This is then used to study infinitesimal deformation of compatible Hom-Leibniz algebra. Furthermore, we establish the relationship between the Nijenhuis operator and the trivial infinitesimal deformation. Further, we introduce the cohomology of compatible Hom-Leibniz algebra with coefficients in an arbitrary representation. 

This paper is organized as follows:
In section $2$,  we recall  some related  basic concepts such as  Hom-Leibniz algebra,  representation of Hom-Leibniz algebra  Balavoine bracket, differential graded Lie algebra and characterisation of Hom-Leibniz algebras. In section $3$,  we define compatible Hom-Leibniz algebra and compatible Hom-bimodules. We then construct the graded Lie algebra whose Maurer-Cartan elements characterize compatible Hom-Leibniz algebra structure.  In section $4$, firstly we introduce  cohomology of compatible Leibniz algebra with coefficients in itself. Then  the cohomology of compatible Leibniz algebra with coefficients in an arbitrary representation  compatible Hom-bimodule) is introduced.
In section $5$,  infinitesimal deformation of compatible Hom-Leibniz  algebra is studied using cohomology of compatible Hom-Leibniz algebra with coefficients in itself. It is shown that equivalent infinitesimal deformations are in the same cohomology group. Then the notion of the Nijenhuis operator on a compatible Hom-Leibniz algebra is studied and the correspondence between the Nijenhuis operator and a trivial deformation is established.  

Throughout the paper we consider the underlying field $K$ to be of characteristic $0$.
\section{Hom-Leibniz Algebra and its characterisation}
In this section, we review the basics of Hom-Leibniz algebras and the Balavoine bracket. Our main references are \cite{nhomLeib, Bala, review}.
\begin{defn} A Hom-Leibniz algebra is a vector space $L$ together with linear operations $[.,.]: L\otimes L \to L$ and  $\alpha: L \to L$ such that
	$$[\alpha(x),[y,z]]=[[x,y],\alpha(z)]+[\alpha(y),[x,z]], ~\forall x,y, z \in L.$$
 \end{defn}
A Hom-Leibniz algebra given by the triple $(L,[~~],\alpha)$ is called multiplicative if $\alpha([x,y])=[\alpha(x),\alpha(y)]$. Hereon  we consider our Hom-leibniz algebras to be multiplicative.
	\begin{defn} A homomorphism between two Hom-Leibniz algebras $(L_1,[~]_1,\alpha_1)$ and $(L_2,[~]_2,\alpha_2)$ is a $K$-linear map $\phi:L_1\to L_2$ satisfying
$$\phi([x,y]_1)=[\phi(x),\phi(y)]_2 ~~\text{and}~~ \phi\circ\alpha_1=\alpha_2\circ \phi.$$
\end{defn}
\begin{defn} Let $(L,[~],\alpha)$ be a Hom-Leibniz algebra. A $L$-bimodule or representation is a vector space $M$ together with  $L$-actions
$m_L:L\otimes M\to M,~~ m_R: M\otimes L\to M$ and a map $\beta\in End(M)$
such that for any $x,y\in L$ and $m\in M$ we have
$$\beta(m_L(x,m))=m_L(\alpha(x),\beta(m)),\;\; \beta(m_R(m,x))=m_R(\beta(m),\alpha(x))$$
$$m_L(\alpha(x),m_L(y,m))=m_L([x,y],\beta(m))+m_L(\alpha(y),m_L(x,m))$$
$$m_L(\alpha(x),m_R(m,y))=m_R(m_L(x,m),\alpha(y))+m_R(\beta(m),[x,y])$$	
$$m_R(\beta(m), [x,y])=m_R(m_R(m,x),\alpha(y))+m_L(\alpha(x), m_R(m,y)).$$
\end{defn}

The following is a well known result.

\begin{prop} Let $(L,[~],\alpha)$ be a Hom-Leibniz algebra and $(M,m_L,m_R,\beta)$ its representation. Then $L\oplus M$ is a Hom-Leibniz algebra with the linear homomorphism $\alpha\oplus \beta: L\oplus M \to L\oplus M$ defined as $(\alpha\oplus \beta)(x,m)=(\alpha(x),\beta(m))$ and the  Hom-Leibniz bracket defined as
$$[(x,u),(y, v)]_\ltimes=([x,y],m_L^1(x,v)+m_R^1(u,y))~~\forall~~x,y\in L~\text{and}~u,v\in M.$$
This is known as the semi-direct product of $L$ and $M$.
\end{prop}
	\begin{defn} A permutation $\sigma \in S_n$ is called an $(i,n-i)$-shuffle if $\sigma(1)<\sigma(2)<...<\sigma(i)$ and $\sigma(i+1)<\sigma(i+2)<...<\sigma(n)$. If $i=0 ~\text{or}~ n$, we assume $\sigma=id$. $S_{(i, n-i)}$ denotes the set of all $(i,n-i)$-shuffles.
	\end{defn}
	\begin{defn} Let $(\mathfrak g=\oplus _{k\in \mathbb Z}\mathfrak g^k,[~], d)$ be a differential graded Lie algebra. A degree 1 element $x\in \mathfrak g^1$ is called a Maurer-Cartan element of $\mathfrak g$ if it satisfies
	$$dx+\frac{1}{2}[x,x]=0.$$
	\end{defn}
	\begin{thm} \cite{review} Let $(\mathfrak g=\oplus _{k\in \mathbb Z}\mathfrak g^k,[~])$ be a graded Lie algebra and $\mu \in \mathfrak{g^1}$ be a Maurer-Cartan element. Then the map
	$$d_\mu:\mathfrak g\to \mathfrak g,~d_{\mu}(u):=[\mu,u], ~\forall u \in \mathfrak g,$$
	 is a differential on $\mathfrak g$. \\
	 Further, for any $v \in \mathbb \mathfrak g^1$, the sum $\mu+v$ is a Maurer-Cartan element of the graded Lie algebra $(\mathfrak g=\oplus _{k\in \mathbb Z}\mathfrak g^k,[~])$ iff $v$ is a Maurer-Cartan element of the differential graded Lie algebra $(\mathfrak g=\oplus _{k\in \mathbb Z}\mathfrak g^k,[~], d_\mu).$
\end{thm}

Let $\mathfrak {g}$ be a vector space and $\alpha:g\to g$ a linear map. For each $n\geq 1$, we denote $ \mathbb C_\alpha^n(\mathfrak {g}, \mathfrak {g})=\{f\in Hom(\otimes ^n\mathfrak {g}, \mathfrak {g}) | \alpha\circ f= f\circ  \alpha^{\otimes n}\}$ and set $\mathbb{C}_\alpha^*(\mathfrak{g}, \mathfrak {g})=\oplus_{n\in \mathbb N}\mathbb C_\alpha ^n(\mathfrak {g}, \mathfrak {g})$.   
	
	We assume the degree of an element in $\mathbb C_\alpha^n(\mathfrak {g}, \mathfrak {g})$ is $n-1$.
	
	For $P\in \mathbb C_\alpha^{p+1}(\mathfrak g, \mathfrak g), Q\in \mathbb C_\alpha^{q+1}(\mathfrak g, \mathfrak g)$ we define the \textbf{ Balavoine bracket} as
	$$[P,Q]_B=P\circ Q-(-1)^{pq}Q\circ P$$
	where $P \circ Q \in \mathbb C_\alpha^{p+q+1} $ is defined as \\
	$$(P\circ Q)(x_1,x_2,...,x_{p+q+1})=\sum_{k=1}^{p+1}(-1)^{(k-1)q}P\circ_k Q,$$
	and 
	$$	Po_kQ(x_1,x_2,...,x_{p+q+1})$$
	{\scriptsize $$	=\sum_{\sigma \in S(k-1, q)}(-1)^{\sigma}P(\alpha^p(x_{\sigma(1)}),...,\alpha^p(x_{\sigma(k-1)}), Q(x_{\sigma(k)},...,x_{\sigma(k+q-1)},x_{k+q}), \alpha^p(x_{k+q+1}),...,\alpha^p(x_{p+q+1}).$$}
	 
	 \begin{thm} \cite{nhomLeib}\label{I} The graded vector space $\mathbb C_\alpha^*(\mathfrak g, \mathfrak g)$ equipped with the Balavoine bracket given above is a graded Lie algebra.
\end{thm}
 
In particular for $\pi \in \mathbb C_\alpha^2(\mathfrak g, \mathfrak g)$, we have $[\pi,\pi]_B\in \mathbb C_\alpha^3(\mathfrak g, \mathfrak g) $ such that\\
$[\pi,\pi]_B=\pi \circ \pi -(-1)^{1.1}\pi \circ \pi= 2 \pi \circ \pi= 2\sum_{k=1}^{2}(-1)^{k-1}\pi \circ_k \pi= 2(\pi \circ_1 \pi - \pi \circ _2 \pi)$ \\
$\pi \circ_1 \pi(x,y,z)=\pi(\pi(x,y),\alpha(z))$ and $\pi \circ_2 \pi(x,y,z)=\pi(\alpha(x), \pi(y,z))-\pi(\alpha(y), \pi(x,z))$.

Thus we have the following corollary. 

\begin{cor} \label{MC} $\pi$ defines a Hom-Leibniz algebra structure on $\mathfrak g$ iff $\pi$ is a Maurer-Cartan element of the graded Lie algebra $(\mathbb C_\alpha^*(\mathfrak g, \mathfrak g),[~]_B).$
\end{cor}
\begin{thm} \label{II} Let $(\mathfrak g, \pi,\alpha)$ be a Hom-Leibniz algebra. Then $(\mathbb C_\alpha ^*(\mathfrak g, \mathfrak g), [~], d_\pi)$ becomes a differential graded Lie algebra (dgLa), where $d_\pi:=[\pi,.]_B$.\\
Further, given $\pi' \in \mathbb C_\alpha^2(\mathfrak g, \mathfrak g)$, $\pi+\pi'$ defines a Leibniz algebra structure on $\mathfrak g$ iff $\pi'$ is a Maurer-Cartan element of the dgLa $(\mathbb C_\alpha^*(\mathfrak g, \mathfrak g), [~], d_\pi)$.
\end{thm}
\section{Compatible Hom-Leibniz algebras and its characterisation}
In this section, we define compatible Hom-Leibniz algebras and then define compatible bimodules over them. We then construct a bidifferential graded Lie algebra whose Maurer Cartan elements are compatible Hom- Leibniz algebras.
\begin{defn} A compatible Hom-Leibniz algebra is a quadruple $(L,[~],\{~\},\alpha)$, where $(L,[~],\alpha)$ and $(L,\{~\},\alpha)$ are Hom-Leibniz algebras such that 	$ ~\forall \,x,y, z \in L$ \begin{equation}
[\alpha(x),\{y,z\}]+\{\alpha(x),[y,z]\}=[\{x,y\},\alpha(z)]+\{[x,y],\alpha(z)\}+[\alpha(y),\{x,z\}]+\{\alpha(y),[x,z]\}. \label{CLA}
\end{equation}
\end{defn}
\begin{prop} \label{equi} A quadruple $(L,[~],\{~\},\alpha)$ is a compatible Hom-Leibniz algebra iff $(L,[~],\alpha)$ and $(L,\{~\},\alpha)$ are Hom-Leibniz algebras such that, for any $k_1,k_2$ in $K$, the bilinear operation 
$$\llbracket x,y\rrbracket=k_1[x,y]+k_2\{x,y\}, ~\forall x,y \in L $$
together with the $k$-linear map $\alpha:L\to L$ defines a Hom-Leibniz algebra structure on L.
\end{prop}
\begin{proof}
Let $(L,[~~],\{~~\},\alpha)$ be a compatible Hom-Leibniz algebra. Then by definition itself $(L,[.,.],\alpha)$ and $(L,\{.,.\},\alpha)$ are Hom-Leibniz algebras. 
Further,
\begin{align*}
& \llbracket \llbracket x,y\rrbracket,\alpha(z)\rrbracket+\llbracket \alpha(y),\llbracket x,z\rrbracket \rrbracket =\llbracket k_1[x,y]+k_2\{x,y\},\alpha(z)\rrbracket+\llbracket \alpha(y), k_1[x,z]+k_2\{x,z\} \rrbracket \notag\\
&=k_1[k_1[x,y]+k_2\{x,y\},\alpha(z)]+k_2\{k_1[x,y]+k_2\{x,y\},\alpha(z)\}+\notag\\
&k_1[\alpha(y), k_1[x,z]+k_2\{x,z\}]+k_2\{\alpha(y), k_1[x,z]+k_2\{x,z\}\}\notag \\
&=k_1k_1[[x,y],\alpha(z)]+k_1k_2[\{x,y\},\alpha(z)]+k_2k_1\{[x,y],\alpha(z)\}+k_2k_2\{\{x,y\},\alpha(z)\}+\notag\\
&k_1k_1[\alpha(y),[x,z]]+k_1k_2[\alpha(y),\{x,z\}]+k_2k_1\{\alpha(y),[x,z]\}+k_2k_2\{\alpha(y),\{x,z\}\}\notag\\
&=k_1^2([[x,y],\alpha(z)]+[\alpha(y),[x,z]])+k_2^2(\{\{x,y\},\alpha(z)\}+\{\alpha(y),\{x,z\}\})\notag\\
&k_1k_2([\{x,y\},\alpha(z)]+\{[x,y],\alpha(z)\}+[\alpha(y),\{x,z\}]+\{\alpha(y),[x,z]\})\notag\\
&=k_1^2[\alpha(x),[y,z]]+k_2^2\{\alpha(x),\{y,z\}\}+k_1k_2([\alpha(x),\{y,z\}]+\{\alpha(x),[y,z]\})\notag\\
&=k_1(k_1[\alpha(x),[y,z]]+k_2[\alpha(x),\{y,z\}])+k_2(k_2\{\alpha(x),\{y,z\}\}+k_1\{\alpha(x),[y,z]\})\notag\\
	&=k_1[\alpha(x),k_1[y,z]+k_2\{y,z\}]+k_2\{\alpha(x),k_2\{y,z\}+k_1[y,z]\}\notag\\
	&=k_1[\alpha(x),\llbracket y,z \rrbracket]+k_2\{\alpha(x),\llbracket y,z \rrbracket\}\notag\\
	&= \llbracket \alpha(x),\llbracket y,z\rrbracket \rrbracket.	\notag  
\end{align*}
The converse is straight forward.
\end{proof}
\begin{defn} A homomorphism between two compatible Hom-Leibniz algebras $(L_1,[~]_1, \{~\}_1,\alpha_1)$ and $(L_2,[~]_2, \{~\}_2,\alpha_2)$ is a k-linear map $\phi:L_1\to L_2$ satisfying
$$\phi([x,y]_1)=[\phi(x),\phi(y)]_2~~\phi(\{x,y\}_1)=\{\phi(x),\phi(y)\}_2~~\text{and}~~\phi\circ \alpha_1=\alpha_2\circ \phi. $$
\end{defn}
\begin{defn} Let $(L,[~],\{~\},\alpha)$ be a compatible Hom-Leibniz algebra. A compatible $L$-bimodule is a vector space $M$ together with four $L$-actions
 $$m_L^i:L\otimes M\to M,~~~~ m_R^i: M\otimes L\to M,\;\; i=1,2$$
and a linear map $\beta: M\to M$ such that
\begin{itemize}
\item $(M,m_L^1, m_R^1,\beta)$ is a bimodule over $(L,[~~],\alpha)$.
\item $(M,m_L^2, m_R^2,\beta)$ is a bimodule over $(L,\{~~\},\alpha)$.
\item the following compatibility conditions hold for all $x,y \in L,~m\in M$
\begin{align*}
LLM: &m_L^1(\alpha(x),m_L^2(y,m))+m_L^2(\alpha(x),m_L^1(y,m))=  m_L^1(\{x,y\},\beta(m))+  \\ & m_L^2([x,y],\beta(m))+
   m_L^1(\alpha(y),m_L^2(x,m))+m_L^2(\alpha(y),m_L^1(x,m))\\
LML:&m_L^1(\alpha(x),m_R^2(m,y))+m_L^2(\alpha(x),m_R^1(m,y))=  m_R^1(m_L^2(x,m),\alpha(y))+ \\& m_R^2(m_L^1(x,m),\alpha(y))+
 m_R^1(\beta(m), \{x,y\})+m_R^2(\beta(m), [x,y])\\
MLL: & m_R^1(\beta(m), \{x,y\})+m_R^2(\beta(m), [x,y]) =m_R^1(m_R^2(m,x),\alpha(y))+ \\& m_R^2(m_R^1(m,x),\alpha(y))+  
m_L^1(\alpha(x), m_R^2(m,y))+m_L^2(\alpha(x),m_R^1(m,y))
\end{align*}
\end{itemize}
\end{defn}
We also say that $(M,m_L^1,m_R^1,m_L^2,m_R^2,\beta)$ is a representation of the compatible Hom-Leibniz algebra $(L,[~],\{~\},\alpha)$.\\
\textbf{Note}: Any compatible Hom-Leibniz algebra $(L,[~],\{~\},\alpha)$ is a compatible $L$-bimodule in which $m_L^1=m_R^1=[~]$ and $m_L^2=m_R^2=\{~\}$.\\
The following result can be proved just like the standard case.
\begin{prop} Let $(L,[~],\{~\},\alpha)$ be a compatible Hom-Leibniz algebra and $(M,m_L^1,m_R^1,m_L^2,m_R^2,\beta)$ its representation. Then $L\oplus M$ is a compatible Hom-Leibniz algebra with the linear homomorphism $\alpha\oplus \beta$ and the compatible Hom-Leibniz brackets defined as
$$[(x,u),(y, v)]_\ltimes=([x,y],m_L^1(x,v)+m_R^1(u,y))~~\text{and}$$
$$\{(x,u),(y, v)\}_\ltimes=(\{x,y\},m_L^2(x,v)+m_R^2(u,y))~~\forall~~x,y\in L~\text{and}~u,v\in M.$$
\end{prop}




\begin{defn} \cite{CLieA} Let ($\mathfrak g, [~~],\delta_1)$ and $(\mathfrak g, [~~],\delta_2)$ be two differential graded Lie algebras. We call $(\mathfrak g,[~~], \delta_1, \delta_2)$ a bi-differential graded Lie algebra (b-dgLa) if $\delta_1$ and $\delta_2$ satisfy
$$\delta_1\delta_2+\delta_2\delta_1=0.$$
\end{defn}
It is easy to show the following.
\begin{prop} \cite{CLieA} Let $(\mathfrak g, [~~],\delta_1)$ and $(\mathfrak g, [~~],\delta_2)$ be two differential graded Lie algebras. Then $(\mathfrak g,[~~], \delta_1, \delta_2)$ is a bi-differential graded Lie algebra iff for any $k_1$ and $k_2$ $\in K$, $(\mathfrak g,[~~], \delta_{k_1k_2})$ is a differential graded Lie algebra where $\delta_{k_1k_2}=k_1\delta_1+k_2\delta_2$.
\end{prop}
\begin{defn} Let $(\mathfrak g,[~~], \delta_1, \delta_2)$ be a b-dgLa. A pair $(\pi_1, \pi_2)\in \mathfrak g_1\oplus \mathfrak g_1$ is called a Maurer-Cartan element of the b-dgLa $(\mathfrak g,[~~], \delta_1, \delta_2)$ if $\pi_1$ and $\pi_2$ are Maurer-Cartan elements of the dgLas $(\mathfrak g,[~~], \delta_1)$ and $(\mathfrak g,[~~], \delta_2)$ respectively, and 
$$\delta _2\pi_1+\delta _1\pi_2+[\pi_1,\pi_2]=0. $$
\end{defn}
\begin{prop} A pair $(\pi_1, \pi_2)\in \mathfrak g_1\oplus \mathfrak g_1$ is a Maurer-Cartan element of the b-dgLa $(\mathfrak g,[~~], \delta_1, \delta_2)$ iff for any $k_1,k_2\in K$, $k_1\pi_1+k_2\pi_2$ is a Maurer-Cartan element of the dgLa $(\mathfrak g,[~~], \delta_{k_1k_2})$.
\end{prop}

\begin{thm} \label{MC bdgLa} Let $(L,\alpha)$ be a Hom-vector space and $\pi_1, \pi_2 \in \mathbb C_\alpha ^2(L,L)$. Then $(L, \pi_1, \pi_2,\alpha)$ is a compatible Hom-Leibniz algebra iff $(\pi_1, \pi_2)$ is a Maurer-Cartan element of the b-dgLa $(\mathbb C_\alpha^*(L,L), [~~]_B, \delta_1=0, \delta_2=0)$.
\end{thm}
		\begin{proof}
		 $(L, \pi_1, \pi_2,\alpha)$ is a compatible Hom-Leibniz algebra. gives $(L, \pi_1,\alpha)$ and $(L, \pi_2,\alpha)$ are Hom-Leibniz algebras. Hence we get $[\pi_1, \pi_1]_B=[\pi_2, \pi_2]_B=0$.\\
		Further, $\forall x,y, z \in L$ we have the compatibility condition,
		\begin{align}
			\pi_1(\alpha(x),\pi_2(y,z))+\pi_2(\alpha(x),\pi_1(y,z))&=\pi_1(\pi_2(x,y),\alpha(z))+\pi_2(\pi_1(x,y),\alpha(z))+\notag\\ 
			&\pi_1(\alpha(y),\pi_2(x,z))+\pi_2(\alpha(y),\pi_1(x,z))
			\label{CC}
		\end{align}
		We note that,	
		$[\pi_1, \pi_2]_B=\pi_1\circ \pi_2+\pi_2\circ \pi_1,$ where 
		$\pi_1\circ \pi_2(x,y,z)=(\pi_1\circ_1 \pi_2-\pi_1\circ_2 \pi_2)(x,y,z)=\pi_1(\pi_2(x,y),\alpha(z))-\pi_1(\alpha(x),\pi_2(y,z))+\pi_1(\alpha(y),\pi_2(x,z))$
		and 
		$\pi_2\circ \pi_1(x,y,z)=(\pi_2\circ_1 \pi_1-\pi_1\circ_2 \pi_1)(x,y,z) =\pi_2(\pi_1(x,y),\alpha(z))-\pi_2(\alpha(x),\pi_1(y,z))+\pi_2(\alpha(y),\pi_1(x,z))$. 
		i.e.,
		\begin{align*}
			 [\pi_1, \pi_2]_B(x,y,z)  &=\pi_1(\pi_2(x,y),\alpha(z))-\pi_1(\alpha(x),\pi_2(y,z)) +\pi_1(\alpha(y),\pi_2(x,z))\\
    & +\pi_2(\pi_1(x,y),\alpha(z))-\pi_2(\alpha(x),\pi_1(y,z))+\pi_2(\alpha(y),\pi_1(x,z)).    
		\end{align*}
	Thus we see that $[\pi_1, \pi_2]=0$ is equivalent to the compatibility condition (\ref{CC}).
		
	\end{proof}
	
	\begin{thm} \cite{CLieA} \label {MC I}Let $(\pi_1,\pi_2)$ be a Maurer-Cartan element of the b-dgLa $(\mathfrak g,[~~], \delta_1, \delta_2)$.\\
	Define $d_1:=\delta_1+[\pi_1,\_]$ and $d_2:=\delta_2+[\pi_2,\_]$. Then  $(\mathfrak g,[~~], d_1, d_2)$ is a b-dgLa.\\
	Further, for any $\Tilde{\pi}_1, \Tilde{\pi}_2 \in \mathfrak g_1$, $(\pi_1+\Tilde{\pi}_1, \pi_2+\Tilde{\pi}_2)$ is a Maurer Cartan element of the b-dgLa $(\mathfrak g,[~~], \delta_1, \delta_2)$ iff $(\Tilde{\pi}_1,\Tilde{\pi}_2)$ is a Maurer-Cartan element of the b-dgLa $(\mathfrak g,[~~], d_1, d_2)$.
	\end{thm}
	Let $(L,\pi_1, \pi_2,\alpha)$ be a compatible Hom-Leibniz algebra. From theorems (\ref{MC bdgLa}) and (\ref {MC I}), we conclude the following important results:
	\begin{thm} $(\mathbb C_\alpha^*(L,L),[~~], d_1,d_2)$ is a b-dgLa where  $d_1:=[\pi_1,\_]_B$ and $d_2:=[\pi_2,\_]_B$.
 \end{thm}
	\begin{thm}  For any $\Tilde{\pi}_1, \Tilde{\pi}_2 \in \mathbb C_\alpha^2(L,L)$, $(L,\pi_1+\Tilde{\pi}_1, \pi_2+\Tilde{\pi}_2)$ is a compatible Hom-Leibniz algebra iff $(\pi_1+\Tilde{\pi}_1, \pi_2+\Tilde{\pi}_2)$ is a Maurer Cartan element of the b-dgLa $(\mathbb C_\alpha^*(L,L),[~~]_B, d_1, d_2)$.
	\end{thm}
\section{Cohomology of compatible Hom-Leibniz algebras }
	Let $(L,[~~],\{~~\},\alpha)$ be a compatible Hom-Leibniz algebra with $\pi_1(x,y)=[x,y]$ and $\pi_2(x,y)=\{x,y\}$.
	By theorem (\ref{MC bdgLa}), $(\pi_1,\pi_2)$ is a Maurer-Cartan element of the b-dgLa $(\mathbb C_\alpha ^*(L,L), [~~]_B,0,0)$.\\
	We define the $n$-cochains for $n\geq 1$ as 
	$$LC_\alpha^n(L,L):=\underbrace{\mathbb C_\alpha^n(L,L) \oplus \mathbb C_\alpha^n(L,L)\oplus\cdots\oplus \mathbb C_\alpha^n(L,L)}_{n \;\text{times}}$$ 
	and $d^n:LC_\alpha^n(L,L) \to LC_\alpha^{n+1}(L,L)$ by
	$d^1f=([\pi_1,f]_B,[\pi_2,f]_B)$; for $n\ge 2$ 
	$d^n(f_1,f_2,..., f_n)=(-1)^{n-1}(g_1,\ldots,g_{n+1}),$
	where $g_1=[\pi_1, f_1]_B$, $g_{n+1}=[\pi_2, f_n]_B $,
	$g_i=[\pi_2, f_{i-1}]_B+[\pi_1, f_i]_B$, for $2\le i\le n$, 
	for all $ (f_1,f_2,...,f_n)\in LC_\alpha^n(L,L)$.
	We have  the following theorem. 
	\begin{thm} \label{d2} We have $d^{n+1}\circ d^n=0$.
 \end{thm}
	\begin{proof}
		We first note that since $(\pi_1,\pi_2)$ is a Maurer-Cartan element of the b-dgLa $(\mathbb C_\alpha^*(L,L), [~~]_B,0,0)$ we have
		$[\pi_1,\pi_1]=0,~[\pi_1,\pi_2]=0,~[\pi_2,\pi_2]=0$.\\
		For any $(f_1,f_2,\cdots,f_n)\in LC_\alpha^n(L,L),$ we have 	$$d^{n+1}d^n(f_1,f_2,\cdots,f_n)
			=(-1)^{n-1}d^{n+1} (g_1,\ldots,g_{n+1}),$$
	 where $g_1=[\pi_1, f_1]_B,\; g_{n+1}=[\pi_2, f_n]_B,\;
	g_i=[\pi_2, f_{i-1}]_B+[\pi_1, f_i]_B,$ for $2\le i\le n$
So, we have 		
		$d^{n+1}d^n(f_1,f_2,\cdots,f_n)=-(h_1,\ldots, h_{n+2})$, where
			$$h_1=[\pi_1,[\pi_1,f_1]_B]_B,\;  h_2=[\pi_2,[\pi_1,f_1]_B]_B+[\pi_1,[\pi_2,f_1]_B]_B+[\pi_1,[\pi_1,f_2]_B]_B, $$ 
			$$h_i=[\pi_2,[\pi_2,f_{i-2}]_B]_B+[\pi_2,[\pi_1,f_{i-1}]_B]_B+[\pi_1,[\pi_2,f_{i-1}]_B]_B+[\pi_1,[\pi_1,f_i]_B]_B,$$ $\forall 3\le i\le n$,
			$h_{n+1}= [\pi_2,[\pi_2,f_{n-1}]_B]_B+[\pi_2,[\pi_1,f_{n}]_B]_B)+[\pi_1,[\pi_2,f_{n}]_B]_B$ and  $\;h_{n+2}= [\pi_2,[\pi_2,f_{n}]_B]_B)$. 
			
			Now, $$h_1=\frac{1}{2}[[\pi_1,\pi_1]_B,f_{1}]_B=0,\; h_2=[[\pi_1,\pi_2]_B,f_{1}]_B+\frac{1}{2}[[\pi_1,\pi_1]_B,f_{2}]_B=0,$$
			$$h_i=\frac{1}{2}[[\pi_2,\pi_2]_B,f_{i-2}]_B+[[\pi_1,\pi_2]_B,f_{i-1}]_B+\frac{1}{2}[[\pi_1,\pi_1]_B,f_{i}]_B=0,\; \forall 3\le i\le n$$
			$$h_{n+1}= \frac{1}{2}[[\pi_2,\pi_2]_B,f_{n-1}]_B+[[\pi_1,\pi_2]_B,f_{n}]_B=0, \;h_{n+2}=\frac{1}{2}[[\pi_2,\pi_2]_B,f_{n}]_B)=0.$$
Thus 	$d^{n+1}\circ d^n=0$.	
		\end{proof}
		Define  $LC_\alpha^*(L,L)=\oplus_{n\in \mathbb N}LC_\alpha^n(L,L)$ and  $d^*=\oplus_{n\in \mathbb N} d^n.$ From above theorem,  $(LC_\alpha^*(L,L), d^*)$ is a cochain complex.
	
	 Let $(L,[~~], \{~~\},\alpha)$ be a compatible Hom-Leibniz algebra. The cohomology of the cochain complex $(LC_\alpha^*(L,L), d^*)$ is called the cohomology of the compatible Hom-Leibniz algebra $(L,[~~], \{~~\},\alpha)$. We denote the corresponding cohomology group by $H_\alpha^n(L,L).$
	
	

Let $g_1$ and $g_2$ be vector spaces. 
For vector spaces $g_1$ and $g_2$, we define $g^{l,k}$ to be the direct sum of tensor products of $g_1$ and $g_2$, where $g_1$ is repeated $l$ times and $g_2$ is repeated $k$ times. For example $g^{1,1}=(g_1\otimes g_2) \oplus(g_2\otimes g_1) $ and $g^{2,1}=(g_1\otimes g_1\otimes g_2)\oplus(g_1\otimes g_2\otimes g_1)\oplus(g_2\otimes g_1\otimes g_1) $. Thus $\otimes ^n(g_1\oplus g_2)\equiv \oplus_{l+k=n}g^{l,k}$.\\ 
For any linear map $f:g_{i_1}\otimes g_{i_2}\cdots \otimes g_{i_n}\to g_j, \text{ where  } i_1,i_2,...,i_n,j\in \{1,2\}$, we define $\hat{f}\in C^n(g_1\oplus g_2,g_1\oplus g_2 )$ as
\[
    \hat{f}= 
\begin{cases}
    f,& \text{on } g_{i_1}\otimes g_{i_2}\cdots \otimes g_{i_n}\\
    0,              & \text{otherwise.}
\end{cases}
\]
$\hat{f}$ is called a lift of $f$.\\
In particular, for the linear maps we encountered in the previous sections:
$$\pi:L\otimes L\to L,~~m_L:L\otimes M \to M,~~m_R:M\otimes L \to M$$
we get lifts
$$\hat{\pi}:\otimes^2(L\oplus M)\to L\oplus M,~\text{defined as}~ \hat{\pi}((x_1,v_1),(x_2,v_2))=(\pi(x_1,x_2),0) $$
$$\hat{m}_L:\otimes^2(L\oplus M)\to L\oplus M,~\text{defined as} ~\hat{m}_L((x_1,v_1),(x_2,v_2))=(0,m_L(x_1,v_2))$$
$$\hat{m}_R:\otimes^2(L\oplus M)\to L\oplus M,~\text{defined as} ~\hat{m}_R((x_1,v_1),(x_2,v_2))=(0,m_R(v_1,x_2))$$

By property  of the Hom-functor we get
$$C^n(g_1\oplus g_2, g_1\oplus g_2)\equiv \sum_{l+k=n}C^n(g^{l,k},g_1)\oplus \sum_{l+k=n}C^n(g^{l,k},g_2). $$

\begin{defn} A linear map $f\in Hom(\otimes ^n(g_1\oplus g_2), (g_1\oplus g_2) )$ has bidegree $l|k$ if
\begin{enumerate}
\item $l+k+1=n$
\item if $X\in g^{l+1,k}$ then $f(X)\in g_1$
\item if $X\in g^{l,k+1}$ then $f(X)\in g_2$
\item $f(X)=0$ in all other cases.
\end{enumerate}
\end{defn}
We use notation $\Vert f \Vert=l|k$. We say that 
$f$ is homogeneous if $f$ has a bidegree.

Considering examples above, we have $\Vert \hat{\pi}\Vert=\Vert \hat{m}_L\Vert=\Vert \hat{m}_R\Vert=1|0$.\\

We consider a few standard results regarding bidegrees.

\begin{lem} \label{L1} If $f_1,f_2,\cdots f_k \in C^n(g_1\oplus g_2, g_1\oplus g_2)$ be homogeneous linear maps and the bidegrees of $f_i$'s are different. Then $f_1+f_2+...+f_k=0$ iff $f_1=f_2=\cdots=f_k=0$.
\end{lem}
\begin{lem} \label{L2} If $\Vert f\Vert =-1|l ~(l|-1)$ and  $\Vert g\Vert =-1|k ~(k|-1)$ then $[f,g]_B=0$.
\end{lem}
\begin{lem} \label{L3} $f\in C^n(g_1\oplus g_2, g_1\oplus g_2)$ and $g\in C^m(g_1\oplus g_2, g_1\oplus g_2)$ be homogeneous linear maps with bidegrees $l_f|k_f$ and $l_g|k_g$ respectively. Then $[f,g]_B$ is a linear map of bidegree $l_f+l_f|k_f+k_g$.
\end{lem}
The following is the hom version of the corresponding result given in \cite{review}.

\begin{thm}  Let $(L,\pi=[~~],\alpha)$ be a Hom-Leibniz algebra. $(M,m_L, m_R,\beta)$ is a representation of $L$ iff  $\hat{m}_L+\hat{m}_R$ is a Maurer Cartan element of the dgLA $(C_{\alpha\oplus \beta}^*(L\oplus V, L\oplus V),[~~]_B, \partial_{\hat{\pi}}=[\hat{\pi}, .]_B)$.
\end{thm}
\begin{cor} \label{MC1} If $(M,m_L, m_R,\beta)$ is a representation of $(L,\pi,\alpha)$, then $[\hat{\pi}+\hat{m}_L+\hat{m}_R,\hat{\pi}+\hat{m}_L+\hat{m}_R]_B=0$. 
\end{cor}
Let $(M,m_L,m_R,\beta)$ be a representation of the Hom-Leibniz algebra $(L,\pi,\alpha)$.\\
We define the set of $n$-cochains as
$$\mathbb C_{\alpha\oplus \beta}^n(L,V):=C_{\alpha\oplus \beta}^{n|-1}(L\oplus M, L\oplus M)\cong C_{\alpha\oplus \beta}(\otimes ^n L,V)$$
and coboundary operator $d^n_{\pi+m_L+m_R}:\mathbb C_{\alpha\oplus \beta}^n(L,M)\to \mathbb C_{\alpha\oplus \beta}^{n+1}(L,M)$ as
$$d^n_{\pi+m_L+m_R}f:=(-1)^{n-1}[\hat{\pi}+\hat{m}_L+\hat{m}_R,\hat{f}]_B,~\forall f\in \mathbb C_{\alpha\oplus \beta}^n(L,V).$$
Note that since $\hat{\pi}+\hat{m}_L+\hat{m}_R\in C^{1|0}$ and $\hat{f}\in C^{n|-1}$, Lemma (\ref{L3}) gives us that $[\hat{\pi}+\hat{m}_L+\hat{m}_R,\hat{f}]_B\in C^{n+1|-1}$.\\

Further note that
$d^{n+1}d^nf=-[\hat{\pi}+\hat{m}_L+\hat{m}_R,[\hat{\pi}+\hat{m}_L+\hat{m}_R,\hat{f}]_B]_B=0$ by the graded Jacobi identity.\\
Thus we have a well defined cochain complex $(\mathbb C_{\alpha\oplus \beta}^*(L,M), d^*_{\hat{\pi}+\hat{m}_L+\hat{m}_R})$.
~~~~~~~~~~~~~~~~~~~~~~~~~~~~~~~~~~~~~~~~~~~~~~~~~~~~~~~~~~~~~~~~~~~~~~~~~~~~~~~~~~~~~~~~~~~~
\begin{thm} Let $(L,\pi_1=[~~],\pi_2=\{~~\},\alpha)$ be a compatible Hom-Leibniz algebra and $(M,m^1_L,m^1_R,m^2_L,m^2_R,\beta )$ a representation of $L$. Then $(\hat{\pi}_1+\hat{m}^1_L+\hat{m}^1_R, \hat{\pi}_2+\hat{m}^2_L+\hat{m}^2_R)$ is a Maurer Cartan element of the bi-differential graded Lie Algebra $(\mathbb C_{\alpha\oplus \beta}^*(L\oplus M, L\oplus M), [~~]_B, 0, 0)$ i.e 
\begin{equation}\label{I}
[\hat{\pi}_1+\hat{m}^1_L+\hat{m}^1_R, \hat{\pi}_1+\hat{m}^1_L+\hat{m}^1_R]_B=0,
\end{equation}
\begin{equation}\label{2}
[\hat{\pi}_1+\hat{m}^1_L+\hat{m}^1_R, \hat{\pi}_2+\hat{m}^2_L+\hat{m}^2_R]_B=0,
\end{equation}
\begin{equation}\label{3}
[\hat{\pi}_2+\hat{m}^2_L+\hat{m}^2_R, \hat{\pi}_2+\hat{m}^2_L+\hat{m}^2_R]_B=0
\end{equation}
\end{thm}
\begin{proof}
 Since  $(M,m^1_L,m^1_R,\beta)$ is a representation of the Hom-Leibniz algebra $(L,\pi_1,\alpha)$, by corollary (\ref{MC1})
equation (\ref{I}) holds. 
Likewise $(M,m^2_L,m^2_R,\beta)$ is a representation of the Hom-Leibniz algebra $(L,\pi_2,\alpha)$, by corollary (\ref{MC1})
equation (\ref{3}) holds.\\
For $x_1,x_2,x_3 \in L$, $v_1, v_2, v_3 \in V$
\begin{eqnarray}
&&[\hat{\pi_1}+\hat{m}^1_L+\hat{m}^1_R, \hat{\pi_2}+\hat{m}^2_L+\hat{m}^2_R]_B(x_1,v_1),(x_2,v_2),(x_3,v_3)\notag\\
	&=&(\pi_1(\pi_2(x_1,x_2),\alpha(x_3)), m^1_L(\pi_2(x_1,x_2),\beta(u_3))+m^1_R(m^2_L(x_1,u_2)+(m^2_R(u_1,x_2),\alpha(x_3)))\notag\\
 &+&(-\pi_1(\alpha(x_1),\pi_2(x_2,x_3),-m^1_L(\alpha(x_1),m^2_L(x_2,u_3)+m^2_R(u_2,x_3))-m^1_R(\beta(u_1),\pi_2(x_2,x_3)))\notag\\
 &+&(\pi_1(\alpha(x_2),\pi_2(x_1,x_3)),m^1_L(\alpha(x_2),m^2_L(x_1,u_3)+m^2_R(u_1,x_3))+m^1_R(\beta(u_2),\pi_2(x_1,x_3)))\notag\\
  &+&(\pi_2(\pi_1(x_1,x_2),\alpha(x_3)),m^2_L(\pi_1(x_1,x_2),\beta(u_3))+m^2_R(m^1_L(x_1,u_2)+m^1_R(u_1,x_2),\alpha(x_3))))\notag\\
   &+&(-\pi_2(\alpha(x_1),\pi_1(x_2,x_3)),-m^2_L(\alpha(x_1),m^1_L(x_2,u_3)+m^1_R(u_2,x_3))-m^2_R(\beta(u_1),\pi_1(x_2,x_3)))\notag\\
    &+&((\pi_2(\alpha(x_2),\pi_1(x_1,x_3)),m^2_L(\alpha(x_2),m^1_L(x_1,u_3)+m^1_R(u_1,x_3))+m^2_R(\beta(u_2),\pi_1(x_1,x_3)))\notag\\
    &=&0.
\end{eqnarray}
We get the above by the compatibility conditions (\ref{CLA}), $LLM$, $LML$ and $MLL$.
\end{proof}
Note that the  coboundary operator for $(L,\pi_1,\alpha)$ with coefficients in $(L,m^1_L,m^1_R,\beta)$ and for $(L,\pi_2,\alpha)$ with coefficients in $(L,m^2_L,m^2_R,\beta)$ are respectively given by
$$d^n_{\pi^1+m^1_L+m^1_R}f:=(-1)^{n-1}[\hat{\pi}_1+\hat{m}^1_L+\hat{m}^1_R,\hat{f}]_B,~\text{and}$$
$$d^n_{\pi^2+m^2_L+m^2_R}f:=(-1)^{n-1}[\hat{\pi}_2+\hat{m}^2_L+\hat{m}^2_R,\hat{f}]_B,~\forall f\in \mathbb C_{\alpha\oplus \beta}^n(L,L).$$
By the graded Jacobi identity it can be shown that the three  conditions (\ref{I}) , (\ref{2}), (\ref{3}) imply
\begin{eqnarray}
    d^{n+1}_{\hat{\pi}_1+\hat{m}^1_L+\hat{m}^1_R}d^n_{\hat{\pi}_1+\hat{m}^1_L+\hat{m}^1_R}=0\notag\\
    d^{n+1}_{\hat{\pi}_2+\hat{m}^2_L+\hat{m}^2_R}d^n_{\hat{\pi}_2+\hat{m}^2_L+\hat{m}^2_R}=0\notag\\
    d^{n+1}_{\hat{\pi}_1+\hat{m}^1_L+\hat{m}^1_R}d^n_{\hat{\pi}_2+\hat{m}^2_L+\hat{m}^2_R}+d^{n+1}_{\hat{\pi}_2+\hat{m}^2_L+\hat{m}^2_R}d^n_{\hat{\pi}_1+\hat{m}^1_L+\hat{m}^1_R}=0\label{*}
\end{eqnarray}
For $n\geq1 $ we define the space of n-cochains $ LC^n(L,M)$ as
$$LC^n(L,M)=\underbrace{\mathbb C_{\alpha\oplus \beta}^n(L,M)\oplus \mathbb C_{\alpha\oplus \beta}^n(L,M)\oplus \cdots \oplus \mathbb C_{\alpha\oplus \beta}^n(L,M)}_{\text{n-copies} }$$

and coboundary for $n\geq 1$, $\partial^n:LC^n \to LC^{n+1}$ as
$$\partial ^1 f= (d_{\hat{\pi}_1+\hat{m}^1_L+\hat{m}^1_R}f, d_{\hat{\pi}_2+\hat{m}^2_L+\hat{m}^2_R}f),~~\forall f\in \mathbb C_{\alpha\oplus \beta}(L,M)$$
and for $2\leq i\leq n$ and $(f_1,f_2,\cdots, f_n)\in LC^n(L,M)$,
$$\partial ^n(f_1,f_2,\cdots, f_n)=( g_1,\ldots, g_{n+1}),$$  where $g_1=d^n_{\hat{\pi}_1+\hat{m}^1_L+\hat{m}^1_R}f_1$, $g_i= d^n_{\hat{\pi}_2+\hat{m}^2_L+\hat{m}^2_R}f_{i-1}+d^n_{\hat{\pi}_1+\hat{m}^1_L+\hat{m}^1_R}f_i,$ $\forall 2\le i\le n$ and  $g_{n+1}= d_{\hat{\pi}_2+\hat{m}^2_L+\hat{m}^2_R}f_n).$
Using (\ref{*})  it can be shown that $\partial ^2=0$.

Let   $(M,m^1_L,m^1_R, m^2_L,m^2_R,\beta)$ be a representation of a compatible Hom-Leibniz algebra $(L,\pi_1,\pi_2,\alpha)$. The cohomology of the cochain complex $(\oplus_{n=1}^\infty LC^n(L,M),\partial)$ is called the cohomology of $(L,\pi_1,\pi_2,\alpha)$ with coefficient in the representation $(M,m^1_L,m^1_R, m^2_L,m^2_R,\beta)$.
The corresponding $n^{th}$ cohomology group is denoted by $\mathbb H_{\alpha\oplus \beta}^n(L,M)$.

\section{Infinitesimal deformations of compatible Hom-Leibniz algebras}
We now define formal one-parameter deformation of a compatible Hom-Leibniz algebra. Focussing on infinitisimal deformations, we also study some aspects of Nijenhuis operators.
\begin{defn} Let $(L,[~~], \{~~\},\alpha)$ be a compatible Hom-Leibniz algebra.
A \textit{formal one-parameter deformation} of  $L$ is a pair of $k[[t]]$-linear maps
$$\mu_t : L[[t]]\otimes L[[t]]\to L[[t]]~ \text{and}$$
$$m_t : L[[t]]\otimes L[[t]]\to L[[t]] ~\text{such that}:$$
\begin{itemize}
  \item[(a)] 
  $\mu_t(a,b)=\sum_{i=0}^{\infty}\mu_i(a,b) t^i$, ~~~
 $m_t(a,b)=\sum_{i=0}^{\infty}m_i(a,b) t^i$  
 
for all $a,b\in T$,  where $\mu_i, m_i:T\otimes T\to T$ are k-linear and $\mu_0(a,b,c)=[a,b]$ and $m_0(a,b)=\{ab\}.$
 
      \item[(b)] For any $t$, $(L[[t]],\mu_t,m_t,\alpha)$ is a compatible Hom-Leibniz algebra.

\end{itemize}
\end{defn}

\begin{defn} Let $(L,[~~], \{~~\},\alpha)$ be a compatible Hom-Leibniz algebra. Let $\mu_1, m_1\in C_\alpha^2(L,L)$. Define
$$\mu_t(x,y)=[x,y]+t\mu_1(x,y),~~~~m_t(x,y)=[x,y]+t m_1(x,y),~~\forall x,y\in L.$$
If for any t, $(L,\mu_t, m_t,\alpha)$ is a compatible Hom-Leibniz algebra, we say that $(L,\mu_t,m_t,\alpha)$ defines an infinitesimal deformation of $(L,[~~], \{~~\},\alpha)$.\\
We also say that $(\mu_1,m_1)$ generates an infinitesimal deformation of $(L,[~~], \{~~\},\alpha)$.
\end{defn}
For convenience we write $[x,y]=\mu_0(x,y)$ and $\{x,y\}=m_0(x,y)$.\\

By theorem (\ref{MC bdgLa}) we have that $(L,\mu_t, m_t,\alpha)$ is a compatible Hom-Leibniz algebra 
$\iff$
$(\mu_t, m_t)$ is a Maurer Cartan element of $(C_\alpha^*,[~~]_B,0,0)$
$\iff$ \begin{itemize}
	\item $[\mu_t,\mu_t]_B=0$
	\item $[m_t,m_t]_B=0$
	\item $[\mu_t,m_t]_B=0$
	\end{itemize}

$\iff$ \begin{itemize}
	\item $[\mu_0,\mu_0]_B=0,~[\mu_0,\mu_1]_B=0, ~[\mu_1,mu_1]_B=0$
	\item $[m_0,m_0]_B=0,~[m_0,m_1]_B=0,~[m_1,m_1]_B=0$
	\item $[\mu_0,m_0]_B=0,~[\mu_0,m_1]_B+[\mu_1,m_0]_B=0,~[\mu_1,m_1]_B=0$.
\end{itemize}

Reordering the terms and excluding the trivial equations we get that
$(L,\mu_t,m_t,\alpha)$ defines an infinitesimal deformation of $(L,[~~], \{~~\},\alpha)$ iff
$$[\mu_0,\mu_1]_B=0,~~~~[m_0,m_1]_B=0~,~~~~[\mu_0,m_1]_B+[\mu_1,m_0]_B~=0$$
$$[\mu_1,\mu_1]_B=0,~~~~[\mu_1,m_1]_B=0,~~~~[\mu_1,m_1]_B=0~.$$

Note that the first line above implies $d^2(\mu_1,m_1)=0$ i.e $(\mu_1,m_1)$ is a 2-cocycle and the second line implies that $(L,\mu_1,m_1,\alpha)$ is a compatible Hom-Leibniz algebra.

Hence we have the following theorem.
\begin{thm}  Let $(L,[~~], \{~~\},\alpha)$ be a compatible Hom-Leibniz algebra. If $(\mu_1,m_1)\in LC_\alpha^2(L,L)$ generates an infinitesimal deformation then $(\mu_1,m_1)$ is a cocycle.
\end{thm}

\begin{defn} Two infinitesimal deformations $(L,\mu_t,m_t,\alpha)$ and $(L,\mu_t',m_t',\alpha)$ of compatible Hom-Leibniz algebra $(L,[~~],\{~~\},\alpha)$ are said to be equivalent if there exists a linear map $N:L\to L$ such that
$$Id+tN:(L,\mu_t,m_t,\alpha) \to (L,\mu_t',m_t',\alpha)$$
is a compatible Hom-Leibniz algebra homomorphism.
\end{defn}
$Id+tN$ being a compatible Hom-Leibniz algebra homomorphism implies 
\begin{enumerate}
	\item $[x,y]=[x,y]'$
	\item $\mu_1(x,y)-\mu_1'(x,y)=[x,N(y)]+[N(x),y]-N[x,y]$
	\item $N\mu_1(x,y)=\mu_1'(x,N(y))+\mu_1'(N(x),y)+[N(x),N(y)]$
	\item $\mu_1'(N(x),N(y))=0$
	\item $\{x,y\}=\{x,y\}'$
	\item $m_1(x,y)-m_1'(x,y)=\{x,N(y)\}+\{N(x),y\}-N\{x,y\}$
	\item $N m_1(x,y)=m_1'(x,N(y))+m_1'(N(x),y)+\{N(x),N(y)\}$
	\item $m_1'(N(x),N(y))=0$
 \item $N\alpha=\alpha N,~~\forall x,y\in L.$
	\end{enumerate}

$2$ and $6$ gives
\begin{eqnarray*}
	(\mu_1-\mu_1')(x,y)&=&([x,N(y)]+[N(x),y]-N[x,y], \{x,N(y)\}+\{N(x),y\}-N\{x,y\})\notag\\
	&=&([\mu_0,N]_B,[m_0,N]_B)
	=d^1N. \notag
	\end{eqnarray*}
Thus we have the following theorem.
\begin{thm} If two infinitesimal deformations $(L,\mu_t,m_t,\alpha)$ and $(L,\mu_t',m_t',\alpha)$ of a compatible Hom-Leibniz algebra $(L,\mu_0,m_0,\alpha)$ are equivalent then, $(\mu_1,\mu_1')$ and $(m_1,m_1')$ are in the same cohomology class.
\end{thm}
\begin{defn} Let $(L,[~~],\alpha)$ be a Hom-Leibniz algebra. A linear map $N:L\to L$ is said to be a Nijenhuis operator on $L$ if 
$$N([x,N(y)]+[N(x),y]-N[x,y])=[N(x),N(y)]~~\forall x,y \in L$$
and $\alpha N=N\alpha.$
\end{defn}
We define linear $[~~]_N: L\otimes L \to L$ as 
$$[x,y]_N=[x,N(y)]+[N(x),y]-N[x,y].$$
Using the multiplicativity of $\alpha$ and the fact that $N\alpha=\alpha N$, we get $$\alpha[x,y]_N=[\alpha(x),\alpha(y)]_N.$$

$T_{[~~]}N: L\otimes L\to L$ denotes the \textit{Nijenhuis torsion} of $N$ defined as\\
$$T_{[~~]}N(x,y)=N([x,y]_N)-[N(x), N(y)], ~~\forall x,y \in L.$$
When $N$ is a Nijenhuis operator we get that $T_{[~~]}N=0$.
\eg The identity map $I:L\to L$ is a Nijenhuis operator on any Hom-Leibniz algbra $(L,[~~],\alpha)$.

\begin{prop} \label{Nijen} If $N:L\to L$ is a Nijenhuis operator on $(L,[~~],\alpha)$, then $(L,[~~]_N,\alpha)$ is also a Hom-Leibniz algebra. Further $N$ is a Leibniz algebra homomorphism from $(L,[~~]_N,\alpha)$ to $(L,[~~],\alpha)$.
Further $(L,[~~],[~~]_N,\alpha)$ forms a compatible Hom-Leibniz algebras.
\end{prop}
\begin{proof} For every $x,y\in L$ put $[x,y]_N=\pi_N(x,y)$ and $[x,y]=\pi(x,y)$. 

Using the Balavoine bracket we get,
$$[\pi_N,\pi_N]_B(x,y,z)=2(\pi_N(\pi_N(x,y),\alpha(z))-\pi_N(\alpha(x),\pi_N(y,z))+\pi_N(\alpha(y),\pi_N(x,z))).$$
Thus $\pi_N=[~~]_N$ defines a Hom-Leibniz algebra structure on L.

Further, $N([x,y]_N)=[N(x),N(y)]$ and $N\alpha=\alpha N$ follows from the definition of Nijenhuis operator and $[~~]_N$.

To show $(L,[~~],[~~]_N,\alpha)$ is a compatible Hom-Leibniz algebras we first note that $\pi_N=[\pi, N]_B$. For any $k_1$ and $k_2\in K$,
\begin{eqnarray*}
[k_1\pi+k_2\pi_N,k_1\pi+k_2\pi_N]_B&=&k_1k_2([\pi,\pi_N]_B+[\pi_N,\pi]_B)\\
&=&2k_1k_2[\pi,\pi_N]_B\\
&=&2k_1k_2[\pi,[\pi,N]_B]_B\\
&=&0.    
\end{eqnarray*}

\end{proof}

\begin{defn}  $(L,[~~],\{~~\},\alpha)$ be a compatible Hom-Leibniz algebra. A linear map $N:L\to L$ is said to be a Nijenhuis operator on $(L,[~~],\{~~\},\alpha)$ if $N$ is a Nijenhuis operator on the Hom-Leibniz algebras $(L,[~~],\alpha)$ and $(L,\{~~\},\alpha)$.
\end{defn}
\begin{prop} Let $(L,[~~],\{~~\},\alpha)$ be a compatible Hom-Leibniz algebra. The linear map $N:L\to L$ is a Nijenhuis operator on $(L,[~~],\{~~\},\alpha)$ iff for any $k_1,k_2$ in $K$, $N$ is a Nijenhuis operator on the Hom-Leibniz algebra $(L,\llbracket~~\rrbracket,\alpha)$, where
$\llbracket x,y \rrbracket=k_1[x,y]+k_2\{x,y\}, ~\forall x,y \in L $.
\end{prop}
\begin{proof} We have  \begin{eqnarray*}
	T_{\llbracket~~\rrbracket}N(x,y)&=&N(\llbracket x,y \rrbracket_N)-\llbracket N(x),(y) \rrbracket\notag\\
	&=&N(k_1[x,y]+k_2{x,y})-k_1[N(x), N(y)]-k_2\{N(x), N(y)\}\notag\\
	&=&k_1(N([x,y]_N)-[N(x), N(y)])+k_2(N(\{x,y\}_N)-\{N(x), N(y)\})\notag\\
	&=&k_1T_{[~~]}N(x,y)+k_2T_{\{~~\}}N(x,y)\notag
\end{eqnarray*}
Hence we have, 
$$T_{\llbracket~~\rrbracket}N=0~~\text{iff}~~~T_{[~~]}N=T_{\{~~\}}N=0.$$
\end{proof}
\begin{prop} \label{Nijen1} Let $(L,[~~],\{~~\},\alpha)$ be a compatible Hom-Leibniz algebra and $N:L\to L$ is a Nijenhuis operator on $(L,[~~],\{~~\},\alpha)$. Then $(L,[~~]_N,\{~~\}_N,\alpha)$ is also a compatible Hom-Leibniz algebra and $N$ is a compatible Hom-Leibniz algebra homomorphism from $(L,[~~]_N,\{~~\}_N,\alpha)$ to $(L,[~~],\{~~\},\alpha)$.
\end{prop}
\begin{proof} Let $N:L\to L$ be a Nijenhuis operator on $(L,[~~],\{~~\},\alpha)$. then by the previous theorem  $N$ is a Nijenhuis operator on the Hom-Leibniz algebra $(L,\llbracket~~\rrbracket,\alpha)$ for any $k_1,k_2$ in $K$.
Using proposition (\ref{Nijen}) we get that $(L,\llbracket~~\rrbracket_N,\alpha)$ is a Hom-Leibniz algebra and $N$ is a Leibniz algebra homomorphism from $(L,\llbracket~~\rrbracket_N,\alpha)$ to $(L,\llbracket~~\rrbracket,\alpha)$.\\
Hence we get that $(L,[~~]_N,\{~~\}_N)$ is a compatible Hom-Leibniz algebra. Further, we also get that
$N$ is a compatible Hom-Leibniz algebra homomorphism from $(L,[~~]_N,\{~~\}_N,\alpha)$ to $(L,[~~],\{~~\},\alpha)$.
\end{proof}

\begin{defn} An infinitesimal deformation $(L,\mu_t,m_t,\alpha)$ of compatible Hom-Leibniz algebra $(L,\mu_0,m_0)$ generated by $(\mu_1,m_1)$ is trivial if there exists linear $N: L\to L$ such that $Id+tN:(L,\mu_t,m_t,\alpha) \to (L,\mu_0,m_0,\alpha)$ is a compatible Hom-Leibniz algebra homomorphism.
\end{defn}
Now $Id+tN$ is a compatible Hom-Leibniz algebra homomorphism iff
\begin{enumerate}
\item $\mu_1(x,y)=[x,N(y)]+[N(x),y]-N[x,y]$
\item $m_1(x,y)=\{x,N(y)\}+\{N(x),y\}-N\{x,y\}$
\item $N\mu_1(x,y)=[N(x),N(y)]$
\item $N m_1(x,y)=\{N(x),N(y)\}$
\item $N\alpha=\alpha N.$
\end{enumerate}
$1,3$ and $5$ gives that $N$ is a Nijenhuis operator on $(L,\mu_0,\alpha)$. $2,4$ and $5$ gives that $N$ is a Nijenhuis operator on $(L,m_0,\alpha)$.

Thus we have the following theorem.

\begin{thm}  \label{NO} A trivial infinitesimal deformation of a compatible Hom-Leibniz algebra gives rise to a Nijenhuis operator.
\end{thm}

\begin{thm} A Nijenhuis operator on a compatible Hom-Leibniz algebra $(L,[~~],\{~~\},\alpha)$ gives rise to a trivial deformation.
\end{thm}
\begin{proof} 
Let $N$ be a  Nijenhuis operator on a compatible Hom-Leibniz algebra $(L,[~~],\{~~\},\alpha)$. Take
$$\mu_1(x,y)=[x,N(y)]+[N(x),y]-N[x,y]$$
$$m_1(x,y)=\{x,N(y)\}+\{N(x),y\}-N\{x,y\}$$
for any $x,y\in L$. Then
\begin{eqnarray*}
d^1N(x,y)&=&([\mu_0,N]_B,[m_0,N]_B)(x,y)\\
&=&([x,N(y)]+[N(x),y]-N[x,y], \{x,N(y)\}+\{N(x),y\}-N\{x,y\})\\
	&=&(\mu_1(x,y),m_1(x,y)).
\end{eqnarray*}
i.e., $(\mu_1,m_1)$ is a 2-cocycle.
Further, since $N$ is a Nijenhuis operator on $(L,[~~],\{~~\},\alpha)$, and $\mu_1=[~~]_N$ and $m_1=\{~~\}_N$, by proposition (\ref{Nijen1}) we get that $(L, [~~]_N,\{~~\}_N,\alpha)$ is a compatible Hom-Leibniz algebra.\\
These two statements implies that $(\mu_1,m_1)$ give rise to an infinitesimal deformation of $L$.
Showing that the deformation is trivial is straightforward.
\end{proof}

\end{document}